\documentclass[12pt]{article}
\date{}
\author{}
\title{\LARGE\bf On vertex covers,matchings and random trees}

\usepackage{t1enc}
\usepackage[latin1]{inputenc}
\usepackage[final]{graphicx}
\usepackage{amsfonts}

\newtheorem{theorem}{Theorem}

\newtheorem{corollary}[theorem]{Corollary}
\newenvironment{remark}{\vspace{.25cm}\addtocounter{theorem}{1} 
\noindent {\bf Remark {\thetheorem}}.}{\vspace{.25cm}}

\newcommand{\bin}[2]{{ #1 \choose #2 }}

\begin{document}

\maketitle


\centerline{\large St\'ephane Coulomb\footnote[1]{Email:
    coulomb@spht.saclay.cea.fr} and  Michel Bauer\footnote[2]{Email:
    michel.bauer@cea.fr}} 

\vspace{.3cm}

\centerline{\large Service de Physique Th\'eorique de
  Saclay\footnote[3]{\textit{ Laboratoire de la Direction des Sciences
      de la Mati\`ere du Commisariat \`a l'Energie Atomique, URA2306
      du CNRS}}}

\vspace{.3cm}

\centerline{\large CE Saclay, 91191 Gif sur Yvette, France}

\vspace{.3cm}

\begin{abstract} We study minimal vertex covers and maximal matchings 
  on trees. We pay special attention to the corresponding
  \textit{backbones} i.e. these vertices that are occupied and those
  that are empty in every minimal vertex cover (resp. these egdes that
  are occupied and those that are empty in every maximal matching).
  The key result in our approach is that for trees, the backbones can
  be recovered from a particular tri-coloring which has a simple
  characterization. We give applications to the computation of some
  averages related to the enumeration of minimal vertex covers and
  maximal matchings in the random labeled tree ensemble, both for
  finite size and in the asymptotic regime.
\end{abstract}

\section{Motivations}

For a given (simple : no loops, no multiple edges) graph, finding the
size of a minimal vertex cover or of a maximal matching (see the
beginning of sec.\ref{secMr} for a reminder of definitions), and
counting the number of solutions all fall under the generic name of
combinatorial optimization problems, a field with a long history.

To analyse the average behavior for these questions, the simplest
model is the Erd\"os-Renyi model of random graphs. In this context,
the evaluation of the average size of a maximal matching has been
solved by Karp and Sipser \cite{ks} (see \cite{afp} for refinements
and \cite{bg3} for a physicist approach) in the thermodynamic limit,
i.e.  when both $V$, the vertex set, and ${\mathcal E}$, the edge set
become large, but the ratio $\alpha =2|{\mathcal E}|/|V|$ has a finite
limit.  For the average size of a minimal vertex cover, the answer is
known only when $\alpha \leq e$ and asymptotically for large $\alpha$
\cite{frieze}.

\vspace{.3cm}

To get more detailed informations on these problems, one can
investigate several combinatorial patterns. In this paper, we
concentrate on backbones, i.e. these vertices that are occupied and
those that are empty in every minimal vertex cover (resp. these egdes
that are occupied and those that are empty in every maximal matching).
While for general graphs the relationships between the backbones are
complicated, this is not true for trees : in that case the backbone
geometry can be recovered from a special tricoloring, unique for each
tree and which is easily characterized.  This is the content of
Theorem \ref{thequiv}, the crucial ingredient for our subsequent
analysis.  We use this theorem to compute the average size of the
backbones and the average number of minimal vertex covers and maximal
matchings for random labeled trees of size $n$, where each of the
$n^{n-2}$ labeled trees has the same probability. We also analyse the
asymptotic behavior for large $n$, see Theorems \ref{enum},\ref{enum2}
and \ref{enum3}.

Due to the simple, locally treelike, structure of Erd\"os-Renyi random
graphs\footnote{For $\alpha \leq 1$, an Erd\"os-Renyi random graph is
  a forest for most thermodynamical purposes, but the local treelike
  structure remains even after the birth of the giant component.
  Moreover, the finite components of size $n$ are distributed with the
  uniform measure on labeled trees of size $n$ in the thermodynamic
  limit.}, insight can often be obtained from an analysis of trees,
even if in the case of this paper the extension is nontrivial. The
present study of backbones and the corresponding applications to
random labeled trees can thus be seen as a preliminary step towards
the analysis of their random graph generalizations.

Our motivation to study backbones comes from physics. The adjacency
matrix of a random graph (which is symmetric) can be seen as an
example of a random Hamiltonian, whose average spectrum one would like
to compute. In the infinite $\alpha$ limit, one recovers a semi
circle, but for finite $\alpha$ matters are much more complicated. The
spectrum contains a dense familly of delta peaks plus presumably a
continuous component when $\alpha$ exceeds a threshold. In the case of
the zero eigenvalue, the structure of the eigenvectors can be studied
in detail\footnote{This is indeed an example of a situation where the
  analysis of random trees proved crucial to understand the case of
  the Erd\"os-Renyi model.}. It exhibits interesting phenomena of
localization and delocalization when $\alpha$ varies \cite{bg2}.  We
shall see below that these phase transitions are closely related to
the structure of the backbones. However, their combinatorial
interpretation is still unclear to us.

\vspace{.3cm}

{\bf Aknowledgements} We thank Martin Weigt for an illuminating remark on
backbones that initiated this work.

\section{Main results}

\label{secMr}

A \textit{ vertex cover} of the graph $A=(V,\mathcal E)$ is a subset
of $V$ containing at least one end of each edge in ${\mathcal E}$. We
are interested in minimal vertex covers, i.e. those whose cardinality
is the smallest. The positive (resp.  negative) vertex-backbone of $A$
is the set of vertices which belong to every (resp. no) minimal vertex
cover. The other vertices are called \textit{ degenerate vertices}. An
edge between degenerate vertices is called \textit{ exclusive} if no
minimal vertex cover contains its two extremities.

A \textit{ matching} of $A$ is a set of non-adjacent edges of $A$.
Maximal matchings are those whose cardinality is the largest. The
positive (resp.  negative) edge-backbone is the set of edges which
belong to every (resp. no) maximal matching. The other edges are
called \textit{ degenerate edges}. A vertex for which there is a
maximal matching none of whose edges contains it as an extremity is
called \textit{ optional}. The vertices that are neither optional, nor
an extremity of an edge in the positive backbone are called \textit{
  unavoidable}.

If $A$ is a tree or a forest, one can characterize these objects by
simple properties and compute them recursively. If the $n^{n-2}$
labeled trees on $n$ vertices are chosen at random with the counting
measure, we shall use this to adress questions of the type ``What is
the average size of the edge-backbones ?'' or ``What is the average
number of maximal vertex covers'' in the random labeled tree ensemble.

Note that we are interested in global extrema. A local version for,
say, minimal vertex covers would be vertex covers such that changing
the state of any occupied vertex to the empty state destroys the
vertex cover property\footnote{The example of a starlike tree shows
  the difference between the local and global versions.}. Some
problems analogous to the ones we deal with but for local problems can
be found for instance in the work of Meir and Moon, see e.g.
\cite{mm} and references therein. For local extrema problems, the
notion of backbones seems to be less relevant.

A \textit{ tricoloring} of the graph $A=(V,\mathcal E)$ is a triple
$(B,{\mathcal R},G) \subset V \times {\mathcal E} \times V$, such that
$B$,$G$ and the set of end-vertices of ${\mathcal R}$ form a partition
of $V$. As a starting point,
\begin{theorem}
\label{thequiv} Suppose $A$ is a tree. Each of the three properties
\textit{(i)},\textit{(ii)} and \textit{(iii)} characterizes one and
the same tricoloring $(B,{\mathcal R},G)$ of $A$.

\textit{(i) Minimal vertex-covers :} $B$ is the positive backbone;
$\mathcal R$ is the set of exclusive edges; $G$ is the negative
backbone.

\textit{(ii) Maximal matchings :} $B$ is the set of unavoidable
vertices; ${\mathcal R}$ is the positive backbone; $G$ is the set of
optional vertices.

\textit{(iii)} The edges in ${\mathcal R}$ are non-adjacent; the edges
with one end-vertex in $G$ have the other end-vertex in $B$; each
vertex in $B$ is connected to $G$ by at least two edges.
\end{theorem}

This unique tricoloring $(B,{\mathcal R},G)$ is called the \textit{
  b-coloring} of $A$. An edge is said \textit{ red} if it is in
$\mathcal R$, and a vertex is said \textit{ brown} if it lies in
$B$,\textit{ green} if it lies in $G$ and \textit{ red} if it is an
end-vertex of a red edge.

In the sequel, $N_c(A)$ denotes the number of vertices with color $c$
(where $c$ is either brown, red or green) in the tree $A$ and $N_c(n)$
is the total number of vertices with this same color among the
$n^{n-2}$ labeled trees on $n$ vertices. We shall work with the
generating functions $F_c(x) \equiv \sum_{n\geq 1}
\frac{N_c(n)}{n!}x^n$. They all involve the tree generating function
$T(x)\equiv \sum_{n\geq 1} \frac{n^{n-1}}{n!}x^n$. Our main
combinatorial and probabilistic results are contained in the following
theorems.

\begin{theorem}
\label{enum}
The generating functions for the total number of brown, red and green
vertices are
$$
F_B=T(x)+T(-T(x))-T(-T(x))^2 \ \ ;\ F_R = T(-T(x))^2 \ \ ;\ F_G=
-T(-T(x)) $$
and the corresponding explicit first terms, closed
formul{\ae}, and asymptotics
\begin{eqnarray*}
N_B & = & (0,0,3,4,185,1026,30457,362664,10245825,195060070,\cdots) \\
\frac{N_B(n)}{n^{n-1}} & = & 1+\sum_{l=1}^n
\left(\frac{-l}{n} \right)^l \left( \frac{2}{l}-1 \right)\bin{n}{l}
\sim 0.2276096757\cdots \\
N_R & = & (0, 2, 0, 48, 120, 4560, 35700, 1048992, 15514128,
456726240,\cdots) \\ 
\frac{N_R(n)}{n^{n-1}} & = & -2 \sum_{l=1}^n
\left(\frac{-l}{n} \right)^l \left( \frac{1}{l}-1 \right)\bin{n}{l}
\sim 0.4104940676\cdots \\
N_G & = & (1, 0, 6, 12, 320, 2190, 51492, 685496, 17286768, 348213690,\cdots) \\
\frac{N_G(n)}{n^{n-1}} & = & -\sum_{l=1}^n \left(\frac{-l}{n} \right)^l
\bin{n}{l} \sim 0.3618962567\cdots
\end{eqnarray*}
\end{theorem}

\begin{corollary}
\label{cor}
The size (that is, the cardinality) of the minimal vertex covers and
maximal matchings of a tree $A$ is $N_B(A)+N_R(A)/2$, hence the
average fraction of vertices in a vertex cover of a tree on $n$
vertices is $n^{1-n}(N_B(n)+N_R(n)/2) \sim 0.4328567095 $ for large
$n$.
\end{corollary}

Let $N_{vc}(n)$ and $N_m(n)$ denote the total numbers of minimal
vertex covers and of maximal matchings among labeled trees on $n$
vertices. The corresponding generating functions, $F_{vc}(x)\equiv
\sum_{n\geq 1} \frac{N_{vc}(n)}{n!}x^n$ and $F_{m}(x)\equiv
\sum_{n\geq 1} \frac{N_{m}(n)}{n!}x^n$, verify

\begin{theorem}
\label{enum2}
The generating function for the total number of minimal vertex covers
is
$$F_{vc}(x)=(1-U)xe^U-UT(x^2e^{2U})+U-\frac{1}{2}U^2,$$
where $xUe^U=T(x^2e^{2U})(e^{xe^U}-1)$.

$N_{vc}=(1, 2, 3, 40, 185, 3936, 35917, 978160, 14301513,
464105440,\cdots)$
\end{theorem}

\begin{theorem}
\label{enum3}
The generating function for the number of maximal matchings is
$$F_m(x)=-\frac{1}{2}(xe^U+U)^2+(1+Uxe^U)xe^U+U-U^2,$$
where
$U=x^2e^{-x^2e^{2U}+xe^U+3U}$.

$N_m=(1, 1, 6, 24, 320, 3270, 55482, 999656, 21718440,
544829130,\cdots)$
\end{theorem}

\begin{remark} 
  Sketch of the relation with the kernel of the adjacency matrix (see
  \cite{bg2} for details). 
  
  The kernel of the adjacency matrix of a tree is directly related to
  the b-coloring.  First one shows, for instance by induction on the
  size of the tree, that the kernel of $A$ has dimension
  $N_G(A)-N_B(A)$.  Second, one shows that the support \footnote{The
    vectors on which the adjacency matrix acts can be interpreted as
    maps from $V$ to the reals and it makes sense to say that a vector
    vanishes on a given vertex.  The support of a vector is the set of
    vertices on which it does not vanish. The support of a familly of
    vectors is the union of the elementary supports.} of the kernel
  consists of the green vertices.
  
  Moreover, the \textit{maximal} subsets $B',G'$ of $V$ such that
  
  \textit{(iii)}' the edges with one end-vertex in $G'$ have the other
  end-vertex in $B'$ and each vertex in $B'$ is connected to $G'$ by
  at least two edges,
  
  coincide with $B$ and $G$ of the b-coloring of $A$. Thus, in the
  case of trees, maximality allows to define the sets $B$ and $G$
  without mentionning ${\mathcal R}$.
  
  Drawing the edges between $B'$ and $G'$ defines a bicolored
  subforest of $A$. But there is a partial converse to these
  constructions : one can show that, for a general graph, a bicolored
  subforest on $B',G'$ satisfying \textit{(iii)}' allows to define a
  $|G'|-|B'|$ dimensional subspace of the kernel with support $G'$.
  For the Erd\"os-Renyi model, the enumeration of the finite maximal
  bicolored subtrees satisfying \textit{(iii)}' accounts for the full
  dimension of the kernel up to an $o(|V|)$ correction for small or
  large $\alpha$, but there is a window of $\alpha$'s for which
  infinite patterns contribute $O(|V|)$ to the dimension of the
  kernel. These are the localization-delocalization transitions
  alluded to before.
  
  These are the results that motivated us to have a closer look at the
  backbones.

\end{remark}

\section{Proof of theorem \ref{thequiv}}

There are many ways to build a proof, depending on personal tastes,
and the choice of the authors has been subject to many fluctuations.
So it is not unlikely that the reader will spare time finding his own
argument instead of understanding the proof that we propose.

Let $A=(V,\mathcal E)$ be a tree. Property $(ii)$ of theorem
\ref{thequiv} obviously characterizes a unique tricoloring of $A$.
Hence, it suffices to establish that
\begin{itemize}
\item Step 1 : $(B,{\mathcal R},G)$ defined by $(i)$ is a tricoloring,
  and it satisfies $(iii)$;
\item Step 2 : $(iii)\Rightarrow (ii)$.
\end{itemize}

If $A=(\{v\},\emptyset)$ is the isolated vertex, any of the three
assertions $(i,ii,iii)$ defines a unique tricoloring
$(\emptyset,\emptyset,\{v\})$ -- $v$ is green -- and we suppose from
now on that $A$ has at least two vertices.

\subsection{Step 1}
Let $(B,{\mathcal R},G)$ be the triple defined by $(i)$ in theorem
\ref{thequiv}. Because $B,G$ and the set of end-vertices of $\mathcal
R$ are mutually disjoint, proving that a degenerate vertex is the end
of an exclusive edge should ensure that $(B,{\mathcal R},G)$ is a
tricoloring of $A$. Then we check that it satisfies $(iii)$.

Deletion of $v\in V$ and its incident edges leaves $p\geq 1$ trees
$A_1,\cdots,A_p$, with $A_i=(V_i,{\mathcal E}_i)$. Denote by $v_i$ the
unique vertex of $A_i$ which is adjacent to $v$ in $A$.

A vertex cover of $A$ obviously induces a vertex cover on each $A_i$.
Conversely, suppose we are given a vertex cover $C_i$ on each $A_i$,
and denote by $C$ the union of the $C_i$'s. An edge of $A$ is either
in some ${\mathcal E}_i$, in which case it has one end in $C_i$ hence
in $C$, or one of the $p$ edges between $v$ and some $v_i$. Hence
$C\cup \{v\}$ is a vertex cover of $A$, but $C$ is not unless it
contains each of the $v_i$'s.

Now, let us write $n_i$ (resp. $\bar{n}_i$) for the minimal
cardinality of vertex covers of $A_i$ containing (resp. not
containing) $v_i$. As a consequence of the previous remarks, a subset
$C$ of $V$ is a minimal vertex cover of $A$ if and only if one of two
exclusive assertions holds :
\begin{itemize}
\item \textit{ Assertion 1} : $1+\sum_i min(n_i,\bar{n}_i) \leq \sum_i
  n_i$, $C$ contains $v$ and $C$ induces a minimal cover on each
  $A_i$.
\item \textit{ Assertion 2} : $1+\sum_i min(n_i,\bar{n}_i) \geq \sum_i
  n_i$, $C$ does not contain $v$ and $C$ induces on each $A_i$ a
  vertex cover of cardinality $n_i$ containing $v_i$.
\end{itemize}

This gives us constraints for $v$ being or not in some backbone, which
are very informative if we note that $n_i\leq \bar{n}_i+1$ for all
$i$.

Suppose first that $v$ is degenerate. Then $1+\sum_i
min(n_i,\bar{n}_i) = \sum_i n_i$, and this implies that
$n_{i_0}=\bar{n}_{i_0}+1$ for a unique $i_0$. There exists a minimal
vertex cover of $A$ containing $v$, which induces a minimal vertex
cover on $A_{i_0}$, hence does not contain $v_{i_0}$. There exists
also a minimal vertex cover not containing $v$, which obviously
contains $v_{i_0}$ : as was to be proved, $v_{i_0}$ is degenerate, and
$v$ is the end of an exclusive edge $\{v,v_{i_0}\}$. $(B,{\mathcal
  R},G)$ is thus a tricoloring. Now, given $i \neq i_0$, we can find a
minimal vertex cover of $A_i$ containing $v_i$ because $n_i \leq
\bar{n}_i$ and then extend it into a minimal vertex cover of $A$
containing both $v$ and $v_i$. Hence $v$ is actually the end of a
unique exclusive edge, proving that the edges of $\mathcal R$ are not
adjacent.

If $v$ is in the negative backbone, a minimal vertex cover of $A$ does
not contain $v$, hence it contains each $v_i$. So the neighbors of
vertices in $G$ are in $B$.

If $v$ is in the positive backbone, $1+\sum_i min(n_i,\bar{n}_i) <
\sum_i n_i$, which proves the existence of at least two distinct $i$'s
such that $n_i=\bar{n}_i+1$. The corresponding $A_i$'s do not admit
minimal covers containing $v_i$. Since a minimal cover of $A$ contains
$v$, it induces a minimal cover on these $A_i$'s : it does not contain
the (at least two) corresponding $v_i$'s. Hence every vertex in $B$
has at least two neighbors in $G$.

This proves that $(B,{\mathcal R},G)$ in $(i)$ is a tricoloring which
satisfies $(iii)$ and we now come to the proof of $(iii)\Rightarrow
(ii)$.

\subsection{Step 2}
Let $(B,{\mathcal R},G)$ be a tricoloring of $A$ satisfying $(iii)$.
Let $R$ denote the set of end-vertices of edges in $\mathcal R$.

If $B=G=\emptyset$, then $\mathcal R$ is a perfect matching of the
tree $A$. Because a tree admits at most one perfect matching
\footnote{This is clear by induction on the size of the tree, if we
  note that an edge ending at a leaf is contained in any perfect
  matching.}, $\mathcal R$ is the unique maximal matching of $A$.

Relax this assumption, suppose $e_0=\{g_0,b_0\}$ is an edge of $A$
($g_0\in G,b_0\in B$) and let ${\mathcal M}$ be a matching of $A$
containing $e_0$. Obviously, there exist paths of the form
$g_0,b_0,\cdots,g_k,b_k$ ($k\geq 0$) such that $\{g_i,b_i\}\in
\mathcal M$, $g_i\in G$ and $b_i \in B$ for all $i$. Take one with
maximal length. Then $b_k$ has at least one neighbor $g_{k+1} \in
G\setminus \{g_k\}$, which is not the end of an edge in $\mathcal M$
(because our path is maximal, and all edges ending at $g_{k+1}$ have
the other end in $B$). Now, we replace in $\mathcal M$ the $k+1$ edges
$\{g_i,b_i\}$ ($0\leq i\leq k$) by the $k+1$ edges $\{b_i,g_{i+1}\}$.
This leads to a matching with same cardinality as $\mathcal M$, but
not containing the edge $\{b_0,g_0\}$.

As a first consequence, there exist maximal maximal matchings not
containing $g\in G$ as an end-vertex (vertices in $G$ are optional).
Now, let $b\in B$ and suppose $\mathcal M$ is a matching not
containing $b$ as an end-vertex. We show that $\mathcal M$ is not
maximal. If some neighbor $v$ of $b$ is not the end-vertex of any edge
in $\mathcal M$ we can append the edge $\{b,v\}$ to $\mathcal M$. On
the other hand, suppose that some neighbor of $b$, $g_0\in G$, is the
end of an edge in $\mathcal M$. Then, we apply the procedure above to
build a matching $\mathcal M'$ of $A$ with same cardinality as that of
$\mathcal M$, not containing $g_0$ as an end-vertex. But these two
matchings coincide except on some path
$g_0,b_0,\cdots,g_k,b_k,g_{k+1}$, which does not contain $b\neq b_0$.
Hence, $\mathcal M'$ does not contain any edge ending at $b$ or $g_0$,
and we can append the edge $\{b,g_0\}$. Hence vertices in $B$ are
unavoidable.

Anticipating the conclusion of this pragraph, let us call \textit{
  forbidden} edges the edges between two vertices in $B$, between a
vertex in $B$ and one in $R$, and those edges not in $\mathcal R$ with
both ends in $R$. Deletion of the forbidden edges leaves some trees,
and $(B,{\mathcal R},G)$ induces on each of these trees a tricoloring
satisfying $(iii)$, of the form $(\emptyset,{\mathcal R}\cap{\mathcal
  E}_i,\emptyset)$ or $(B\cap V_i,\emptyset,G\cap V_i)$. Moreover, a
matching $\mathcal M$ induces a matching on each of these trees and,
if $\mathcal M$ contains $p\geq 1$ forbidden edges, at least $p+1$ of
these induced matchings do not contain some vertex in $B$ as an
end-vertex or some edge in $\mathcal R$. By the preceding remarks,
they are not maximal. Hence, by deleting from $\mathcal M$ these $p$
edges, and by replacing the edges of these $p+1$ matchings by those of
maximal matchings, we obtain a matching of $A$ containing at least one
more edge than $\mathcal M$. So a maximal matching does not contain
any forbidden edge and, as an easy corollary, deletion of these edges
leaves maximal matchings of the resulting trees.

Thus, edges in $\mathcal R$ are in the positive backbone of $A$.
Denote by $(\hat{B},\hat{\mathcal R},\hat{G})$ the tricoloring of $A$
defined by $(ii)$ : we have proved that $B\subset \hat{B},{\mathcal
  R}\subset \hat{\mathcal R},G\subset \hat{G}$. By general properties
of tricolorings, this implies that $(B,{\mathcal
  R},G)=(\hat{B},\hat{\mathcal R},\hat{G})$ and concludes the proof of
theorem \ref{thequiv}.

\paragraph{Remark} Let us consider a minimal vertex cover $C$ of $A$. It
contains all the $N_B(A)$ brown vertices of $A$ and none of the green
vertices. The other $N_R(A)$ vertices are ends of non-adjacent red
edges, and we have seen that $C$ contains exactly one end of each such
edge : hence $C$ contains exactly $N_B(A)+N_R(A)/2$ vertices.

A maximal matching of $A$ contains all the $N_R(A)/2$ red edges and
exactly one edge ending at each brown vertex, the other end being
green (hence not brown). Moreover it does not contain any other edge :
hence a maximal matching of $A$ contains exactly $N_B(A)+N_R(A)/2$
edges. This proves corollary \ref{cor}.

\section{Generating functions}

\subsection{Generating function for b-colorings}

Our purpose in this section is to give an exponential generating
function for the number of labeled trees with given color
distribution :
$$F(g,b,r) \equiv \sum_{n\geq 1}\sum_{A \in {\mathcal A}_n}
\frac{1}{n!} g^{N_G(A)}b^{N_B(A)}r^{N_R(A)}$$

where ${\mathcal A}_n$ is the set of labeled trees on $n$ vertices.
\\
Recalling that a tree has a unique b-coloring, we say that a rooted tree
has color $c$ if its root has color $c$.  Let $G,B,R$ be the
(exponential) generating functions for respectively green, brown, red
rooted trees.

Let $A$ be a rooted tree. Then
\begin{itemize}
\item $A$ is green if, and only if, its root is connected to the root
  of arbitrarily many trees defined as follows : root adjacent to
  arbitrarily many rooted colored trees, with the condition that at
  least \textit{ one} root be green. Let us call \textit{ quasi-brown}
  these trees and denote by $U$ their generating function. Then
\begin{eqnarray}
\label{eqG}
G & = & ge^U \\
\label{eqU}
U & = & b e^{B+R} (e^G-1)
\end{eqnarray}
\item $A$ is brown if, and only if, its root is connected to
  the root of arbitrarily many brown or red rooted trees and to at
  least \textit{ two} green rooted trees, so

\begin{equation}
B=b e^{B+R} (e^G-1-G)
\label{eqB}
\end{equation}

\item Finally, $A$ is red if, and only if, its root is connected to
  arbitrarily many brown or red rooted trees and to exactly one tree
  defined as follows : root adjacent to arbitrarily many red or brown
  rooted trees. Let us call \textit{ quasi-red} these trees and denote
  by $Q$ their generating function. Then
\begin{eqnarray}
\label{eqR}
R & = & rQe^{B+R} \\
\label{eqQ}
Q & = & re^{B+R}
\end{eqnarray}
\end{itemize}

Now, the generating function $F$ for colored trees is the only
function of $g,b,r$ such that $g \frac{\partial F}{\partial g}=G$, $b
\frac{\partial F}{\partial b}=B$, $r \frac{\partial F}{\partial r}=R$
and $F(0,0,0)=0$. The following $F$ indeed satisfies these conditions,
thus the generating function for colored trees is
\begin{equation}
\label{eqF}
F(g,b,r)=-\frac{1}{2} ((B+R)^2+Q^2)-GU+be^{B+R} (e^G-1-G)+ge^U+rQe^{B+R}
\end{equation}

We check that $F(x,x,x)$ gives back the usual generating function for
labeled trees : $F_0 (x)=\sum_{n \geq 1} \frac{n^{n-2}}{n!} x^n$.
Recall that the generating function $T=xF_0'(x)$ for rooted trees
verifies $T(x)=xe^{T(x)}$ for $|x|<1/e$.
\\
Putting $S=B+R$ and taking $g,b,r=x$ in the equations
(\ref{eqB})+(\ref{eqR})-(\ref{eqU}) and (\ref{eqQ})-(\ref{eqG}) yields
\begin{eqnarray*}
S-U & = & (Q-G)xe^S \\
Q-G & = & xe^U (1-e^{S-U})
\end{eqnarray*}
Taking $x\rightarrow 0$, this implies $S=U$ and $G=Q$, so
(\ref{eqU})+(\ref{eqQ}) yields $S+G=xe^{S+G}$. Hence $S+G=S+Q=T(x)$,
and it follows from (\ref{eqQ}) that $Qe^Q=xe^{S+Q}=xe^{T(x)}=T(x)$.
Finally :

\begin{eqnarray*}
S & = & U=T(x)+T(-T(x)) \\
G &= & Q=-T(-T(x)) \\
\end{eqnarray*}

Inject these into $F(x,x,x)$ to get $F=T(x)-\frac{1}{2} T(x)^2$, which
is indeed equal to $F_0(x)$.

\begin{remark} Sketch of the relation with Feynman graph enumeration.
  
  Define $\mathcal{S}(S=B+R,U,G,Q,g,b,r)$ by the right hand-side of
  eq.(\ref{eqF}) but seen as a function of seven independent
  variables. Then the vanishing of the partial derivative of
  $\mathcal{S}$ with respect to $S$ leads to the combination
  eq.(\ref{eqB})+eq.(\ref{eqR}), whereas the vanishing of the partial
  derivatives with respect to $U,G$ and $Q$ leads to eq.(\ref{eqG}),
  eq.(\ref{eqU}) and eq.(\ref{eqQ}) respectively. Thus, $F$ is the
  value of $\mathcal{S}$ at the (unique in the small $g,b,r$
  expansion) extremum in the capital variables. The same kind of
  considerations would apply to all the generating functions in this
  paper
  
  We do not know if there is a simple combinatorial explanation for
  this extremal property, but there is a simple physical
  interpretation that we give in appendix \ref{app} to illustrate how
  two scientific communities deal with the same problem. The reader
  interested in a more thourough study of the combinatorics of Feynman
  graphs can consult e.g. \cite{iz}.

\end{remark}

The generating function for the total number of vertices with a given
color comes from differentiation with respect to the corresponding
variable, followed by the identification $b=g=r=x$ :
\begin{eqnarray*}
\sum_n \frac{N_B(n)}{n!}x^n & = &  T(x)+T(-T(x))-T(-T(x))^2\\
\sum_n \frac{N_G(n)}{n!}x^n & = & -T(-T(x)) \\
\sum_n \frac{N_R(n)}{n!}x^n & = & T(-T(x))^2
\end{eqnarray*}
In order to give explicit formul{\ae} for the average numbers of
vertices of each color, we need to know the term of given degree in
$T(-T(x))$ and $T(-T(x))^2$. Writing these as contours integrals along
a small contour surrounding 0 and changing the integration variable
$x$ into $-te^{t}$ yields
\begin{eqnarray*}
\oint \frac{dx}{x^{n+1}} T(-T(x)) & = & \frac{(-1)^n}{n}\oint
\frac{dt}{t^n} e^{-nt}T'(t) \\
\oint \frac{dx}{x^{n+1}} T(-T(x))^2 & = & 2\frac{(-1)^n}{n} \left( \oint
\frac{dt}{t^n} e^{-nt}T'(t)-\oint \frac{dt}{t^{n+1}} e^{-nt}T(t) \right),
\end{eqnarray*}
from which follow both the closed forms and, by the steepest descent
method, their large-size asymptotics\footnote{In this simple
  situation, we can proceed naively to get the asymptotics. For a more
  rigorous treatment in a similar but slightly more involded context,
  see e.g. \cite{drmota}.}
\begin{eqnarray*}
\frac{N_B(n)}{n^{n-1}} & = & 1+\sum_{l=1}^n \left( \frac{-l}{n} \right)^l
\left( \frac{2}{l}-1 \right) \bin{n}{l} \sim 1+T'(-1)+2T(-1)) \\
\frac{N_G(n)}{n^{n-1}} & = & -\sum_{l=1}^n \left(\frac{-l}{n} \right)^l
\bin{n}{l} \sim T'(-1)\\
\frac{N_R(n)}{n^{n-1}} & = & -2 \sum_{l=1}^n \left(\frac{-l}{n} \right)^l
\left( \frac{1}{l}-1 \right)\bin{n}{l} \sim -2(T'(-1)+T(-1)) 
\end{eqnarray*}

\subsection{Generating function for minimal vertex covers}

Let us define a \textit{ covered tree} to be a pair $(A,C)$, where $A$ is
a rooted tree and $C$ is a minimal vertex cover of $A$. The generating
functions for brown and green covered trees are denoted respectively
by $B$ and $G$ in this section. For red covered trees, it is useful to
make the distinction between the minimal covers which contain the root
and those which do not : let us denote by $R_+,R_-$ the corresponding
generating functions.

Consider a pair $(A,C)$, where $A$ is a rooted tree with root $v$ and
$C$ a subset of the set of vertices of $A$. Then
\begin{itemize}
\item $(A,C)$ is a green covered tree if, and only if, $v\notin C$,
  $v$ is attached to arbitrarily many quasi-brown trees, and $C$
  induces on each of these trees a vertex cover with minimal
  cardinality among those containing the root.
\item $(A,C)$ is a brown covered tree if, and only if, $v\in C$, $v$
  is attached to at least 2 green trees and to arbitrarily many brown
  or red rooted trees, and $C$ induces on each of these trees a
  minimal vertex cover.
\item $(A,C)$ is a red covered tree if, and only if, $v$ is attached
  to exactly one quasi-red tree $A_{i_0}$ and to arbitrarily many
  brown or red trees, and one of two exclusive assertions holds :
  $(1)$ $v\in C$, $v$ induces a minimal cover on each of the attached
  tree; $(2)$ $v\notin C$, $v$ induces on each of the attached trees a
  cover with minimal cardinality among those containing the root.
\end{itemize}

This leads to
\begin{eqnarray*}
B=b(e^G-1-G)e^{B+R_+ +R_-} & & G=ge^U   \\
R_+=rQ_- e^{B+R_+ +R_-} & & R_-=rQ_+ e^{B+R_+}\\
\end{eqnarray*}
where the auxiliary function $U,Q_+,Q_-$ are defined as
$$U\equiv b(e^G-1)e^{B+R_+ +R_-}, \ \ Q_+\equiv re^{B+R_-+R_+}, \ \ 
Q_-\equiv re^{B+R_+}.$$

The generating function for covered rooted trees is $G+B+R_+ +R_-$ and
\begin{eqnarray*}
F_{vc} & = & ge^U+b(e^G-1-G)e^{B+R_+ +R_-}+rQ_- e^{B+R_+ +R_-} +rQ_+
e^{B+R_+} \\  & & \hspace{1cm} -GU-\frac{1}{2}(B^2+R_+^2)-R_+ R_- -B(R_+
+R_-)-Q_+ Q_-
\end{eqnarray*}

turns out to be the only function with correct $b,r,g$ partial
derivatives satisfying $F_{vc}(0,0,0)=0$.

Let us identify $b,r,g=x$. Then $B+R_+=U$, $G=Q_-$ and $R_+ =R_-
\equiv R=x^2 e^{2U+R}=T(x^2e^{2U})$. Hence the closed formula for $U$
is $xUe^U=(e^{xe^U}-1)T(x^2e^{2U})$, and the expression for $F_{vc}$
follows immediately.

\subsection{Generating function for maximal matchings}
We shall skip the details, the crucial points being that
\begin{itemize}
\item A maximal matching of a tree $A$ contains all the red edges and
  exactly one edge ending at each brown vertex, the other end being
  green. It does not contain any other edge.
\item Given an edge $B-G$, there exist maximal matchings which do not
  contain it (because $g$ is optional). There also exist some which do
  contain it. Indeed, let $e=\{b,g\}$ be an edge of $A$ with $b\in
  B,g\in G$. There exists a maximal matching not containing $g$ as an
  end vertex and, as a maximal one, this matching contains an edge
  $e'$ ending at $b$. Just replace $e'$ by $e$.
\end{itemize}

Then the generating functions for matched trees read
\begin{eqnarray*}
G_+=g U_- e^{U_+} & & G_-=ge^{U_+} \\
B=b G_- (e^{G_+ +G_-}-1)e^{B+R} & & R=r Q e^{B+R}
\end{eqnarray*}
where
$$U_+ \equiv b G_-e^{G_+ +G_- +B+R}, \ \ U_-\equiv b(e^{G_+
  +G_-}-1)e^{B+R}, \ \ Q\equiv r e^{B+R}$$
The generating function
writes 
\begin{eqnarray*}F_m & = & g U_- e^{U_+} + ge^{U_+} +b
  G_- (e^{G_+ +G_-}-1)e^{B+R}+r Q e^{B+R}\\ & & \hspace{1cm}
-G_+ U_+ -G_- (U_+ +U_-)-\frac{1}{2}(B+R)^2 -\frac{1}{2}Q^2.
\end{eqnarray*}
For $b=r=g=x$, we find $B+R=U_+\equiv U,Q=G_-$.  Hence,
$G_+=U-x^2e^{2U}$ and $U=x^2e^{3U+xe^U-x^2e^{2U}}$, as was to be
proved.

\begin{remark} The quantities $\frac{1}{n}\log N_{vc}(n)$,
  $\frac{1}{n}\log N_{m}(n)$ given analytically by theorems
  \ref{enum2},\ref{enum3} are difficult to confront to numerical
  simulations sampling the $n^{n-2}$ trees uniformly. They are typical
  examples of non self-averaging quantities. This means basically
  that a small fraction of trees contributes significantly to the
  average although it is unlikely to be ``visited'' in reasonable time
  by a Monte-Carlo algorithm sampling trees uniformly. A simpler task 
is the numerical estimation of $\frac{1}{n}<\log N_{vc}(A)>$
or $\frac{1}{n}<\log N_{m}(A)>$, a self averaging quantity which
answers the question : how many minimal vertex covers or maximal
matchings does a typical tree have. This question will
 be adressed analytically and numerically in a work to come, again
 based on the use of b-colorings \cite{cool}.
\end{remark}

\appendix
\section{A note on Feynman graphs}
\label{app}
In quantum field theory, graph counting occurs for the following
reasons. A physical system is characterized by an action
$\mathcal{S}(T_a,g_i)$ where the $T_a$'s denote dynamical variables
and $g_i$ coupling constants. The index $a$ often runs through a
continuum, but for the present dicussion, we assume it to take a
finite number of values. This is the usual case that there is only a
finite number of coupling constants. The general structure of
$\mathcal{S}$ is $\mathcal{S}=\mathcal{S}_0+\sum_i g_iP_i(T_a)$ where
$\mathcal{S}_0=-\frac{1}{2}\sum_{a,b}C_{ab}T_aT_b$ is a quadratic form
which we assume here to be nondegenerate and the $P_i(T_a)$ are
analytic functions. The quantity to be computed is the free energy
  
\begin{equation}\label{eqQFT} \hbar \log \int \prod_a
    \frac{dT_a}{\sqrt{2\pi\hbar}} \sqrt{\det C}  \exp
  \frac{\mathcal{S}(T_a,g_i)}{\hbar},\end{equation}

where the integration contours and values of the $g_i$'s are choosen
to ensure convergence of the integral, and $\hbar$ is Plank's
constant. Note that when the coupling constants $g_i$ all vanish, this
expression vanishes too. The so-called semi-classical expansion
expresses the free energy of the system as an asymptotic expansion in
powers of $\hbar$, the $\hbar^n$ term being computable by definite
rules from certain non simple connected graphs, the so-called Feynman
graphs, with $n$ independent cycles. To understand the appearance of
graphs, the easy way is to first expand formally the integrand in
(\ref{eqQFT}) in powers of the coupling constants and then the $P_i$'s
in powers of the fields $T_a$. This reduces the integral to
integration of a monomial against a gaussian weight, and the
combinatorics of the result is obtained by repeated integration by
parts, which amounts to pair successively and in all possible ways
all pairs of variables $T_aT_b$ in the monomial and replace them by
$\hbar(C^{-1})_{ab}$. In that way, one can interpret the monomial in
the $P_i$'s as vertices marked with the fields they involve, and
$\hbar(C^{-1})_{ab}$ is the weight for an edge of type $ab$ between
two vertices : all possible graphs that can be built in that way
appear in the formal power series expansion of the integral. Taking
the logarithm to compute the free energy amounts to keep only
connected graphs as usual in combinatorics.

To make this general idea concrete take $\mathcal{S}=-T^2/2+ge^T$. In
that case the quadratic form is associated to the $1\times 1$ identity
matrix and $e^T$ describes vertices of arbitrary degree. The integral
along the real axis with a purely imaginary $g$ makes sense and on can
compute the asymptotic expansion at small $g$. The result is that
$$
\int \frac{dT}{\sqrt{2\pi\hbar}} \exp
\frac{\mathcal{S}(T,g)}{\hbar}\sim \sum_{n \geq 0}
\frac{(g\hbar^{-1})^n}{n!}e^{n^2\hbar/2}$$
On the other hand, if to an arbitrary graph (not necessarily simple)
on $n$ vertices,  described by a symmetric matrix $M=(m_{pq})$, one gives
a weight $$\prod_{p\leq q}\frac{\hbar^{m_{pq}}}{m_{pq}!}
  \prod_{p}\frac{1}{2^{m_{pp}}}$$ 
which is essentially its symmetry factor, the sum over all graphs
reconstructs the factor $e^{n^2\hbar/2}$. 

The classical limit corresponds to keeping only the $\hbar^0$
contribution, i.e.  trees. On the other hand, in the classical limit
the system is decribed by the classical equations of motion, which say
that the action $\mathcal{S}$ is extremal with respect to all field
variations : in the present context, this boils down to the stationary
phase method. And indeed, for the concrete example above, this
extremum condition leads to $T=ge^T$, so that $T$ is the rooted labeled tree
generating function. 

We leave it to the reader to compute the inverse of $C$, to extract
the kinds of vertices and edges that Feynman graphs produce when
$$\mathcal{S}(T,U,G,Q,g,b,r)=-\frac{1}{2} (T^2+Q^2)-GU+be^T
(e^G-1-G)+ge^U+rQe^T $$
and retrieve in that way the conditions \textit{(iii)}.

The main message is that, while for a combinatorist
eqs.(\ref{eqG},\ref{eqU},\ref{eqB},\ref{eqR},\ref{eqQ}) for rooted
trees follow from
routine arguments, for a quantum field theorist it is
$\mathcal{S}(T,U,G,Q,g,b,r)$ which comes immediately to mind to count the
desired unrooted graphs and, at the extremum in $(T,U,G,Q)$, trees.

\end{document}